\documentclass[a4paper,12pt]{amsart}
\usepackage{amsmath,latexsym,amscd}

\title{A sharp lower bound for the canonical volume of 3-folds of general type}
\author{Meng Chen}
\address{School of Mathematical Sciences, Fudan University, Shanghai, 200433,
People's Republic of China}
\email{mchen@fudan.edu.cn}
\thanks{The author is supported by the National Natural Science Foundation of China.}
\newcommand{\bC}{{\mathbb C}}
\newcommand{\bQ}{{\mathbb Q}}
\newcommand{\bbP}{\mathbb P}
\newcommand{\so}{\mathcal O}
\newcommand{\roundup}[1]{\ulcorner{#1}\urcorner}
\newcommand{\rounddown}[1]{\llcorner{#1}\lrcorner}

\newtheorem{thm}{Theorem}[section]
\newtheorem{lem}[thm]{Lemma}
\newtheorem{cor}[thm]{Corollary}

\newtheorem{claim}[thm]{Claim}
\theoremstyle{definition}
\newtheorem{defn}[thm]{Definition}
\newtheorem{setup}[thm]{}
\newtheorem{question}[thm]{Question}
\newtheorem{exmp}[thm]{Example}

\newtheorem{rem}[thm]{Remark}
\theoremstyle{remark}

\begin{document}
\begin{abstract}
Let $V$ be a smooth projective 3-fold of general type. Denote by
$K^3$, a rational number, the self-intersection of the canonical
sheaf of any minimal model of $V$. One defines $K^3$ as the
canonical volume of $V$. Assume $p_g(V)\ge 2$. We show that
$K^3\ge \frac{1}{3}$, which is a sharp lower bound. Then we
classify those $V$ with small volume. We also give some new
examples with $p_g=2$ which have maximal canonical stability
index. Finally, we give an application to 4-folds of general type.
\end{abstract}
\maketitle
\pagestyle{myheadings} \markboth{\hfill M. Chen\hfill}{\hfill A
sharp lower bound for the canonical volume\hfill}

\section{\bf Introduction}
To classify algebraic varieties is one of the main goals of
algebraic geometry. For a long time, we have been interested in
classifying algebraic 3-folds of general type which are,
naturally, quite important objects for birationalists. The book
\cite{C-R} edited by Corti and Reid explains some ways to
understand the explicit structure of algebraic threefolds. One
might have noted, however, that this is a very big topic and that
even there is not any answer to lots of very elementary questions.

Let $V$ be a smooth 3-dimensional projective variety of general
type. According to Mori's Minimal Model Program (see for instance
a sample of references \cite{KMM, K-M, Reid83}), $V$ has at least
one minimal model $X$ which is normal projective with only
${\bQ}$-factorial terminal singularities. Denote by $K^3:=K_X^3$,
which is a positive rational number. Reid (the last section of
\cite{Reid83}) first showed that all minimal models of $V$ are in
fact isomorphic in codimension 1. One also knows the uniqueness of
the Zariski decomposition for canonical bundles $mK_V$ by Kawamata
(see page 355 of \cite{KMM} and Proposition 4 in
\cite{Ka-zariski}). Therefore one can see that the canonical
volume $K^3$ is uniquely determined by the birational equivalence
class of $V$. $K^3$ is an important invariant and the value
$\text{Vol}(V):=K^3$ is referred to as {\it the canonical volume
of} $V$. Obviously $V$ has some other important birational
invariants such as the geometric genus $p_g(V):=\dim H^0(V,
\Omega_V^3)$, the irregularity $q(V):=\dim H^1(V, {\mathcal
O}_V)$, the second irregularity $h^2({\mathcal O}_V):=\dim H^2(V,
{\mathcal O}_V)$. These invariants determine the holomorphic
Euler-Poincar\'e characteristic
$$\chi({\mathcal O}_V):=1-q(V)+h^2({\mathcal O}_V)-p_g(V).$$
A difficulty arises in the study of 3-folds of general type is
that $K^3$ is only a rational number, rather than an integer.
Furthermore it may be, in fact, very small. For example, among
known ones by Fletcher-Reid (page 151 in \cite{C-R}),
$\text{Vol}(V)$ could be as small as $\frac{1}{420}$. Inspired by
a series of lectures of Y.T. Siu, we consider the following very
natural and interesting question.
\begin{question}\label{Q} What is the sharp lower bound of the volume among all
minimal projective 3-folds of general type?
\end{question}

Since several decades, it is a fact that the treatment of
pluricanonical maps is easier when $K^3$ is larger. Morever, the
treatment of pluricanonical maps is much easier when $p_g>0$. This
is due to the fact that under this assumption there are no gaps,
since a pluricanonical system $|(m+1)K|$ always contains the
pluricanonical system $|mK|$ as a subsystem. {}From this point of
view, the lower bound for $K^3$ becomes quite important. After the
first manuscript of this paper was finished, we noted that Hajime
Tsuji, Christopher D. Hacon and James Mckernan \cite{Tsuji, H-M}
announced the existence of an uniform lower bound of $K^3$. In
this paper, we are only able to answer Question \ref{Q} under the
extra assumption $p_g\ge 2$. We will see that the volume strongly
affects the canonical stability index (see \ref{cs} below for the
definition) of the given 3-fold $V$.

\begin{defn}\label{cs}
The {\it canonical stability index} $n_0(Y)$ of any smooth
projective variety $Y$ of general type is defined to be the
minimal integer $n_0$ such that for all $i\ge n_0$ the i-th
pluricanonical map is birational onto its image. For $p_g>0$, it
is simply the smallest integer $n_0$ such that the $n_0$-th
pluricanonical map is birational onto its image.
\end{defn}

\begin{setup}\label{kr}{\bf Upper bound for $n_0(V)$.} Let $V$ be
a smooth 3-dimensional projective variety of general type. In
\cite{IJM}, the following was proved:

(1) $n_0(V)\le 5$ whenever $p_g(V)\ge 4$ and this bound is sharp;

(2) $n_0(V)\le 6$ whenever $p_g(V)=3$ and this bound is sharp;

(3) $n_0(V)\le 8$ whenever $p_g(V)=2$.

Recently I was informed of separately by both Ezio Stagnaro and
Christopher Hacon \cite{H} that the bound of $n_0(V)$ in
\ref{kr}(3) is also sharp. As being stated in Theorem
\ref{main}(2), more examples with $p_g=2$ and with maximal
canonical stability index can be found if the canonical volume
satisfies certain inequality.
\end{setup}

My main results are as follows.

\begin{thm}\label{main}
Let $V$ be a smooth 3-dimensional projective variety of general
type. Assume $p_g(V)\ge 2$. Then

(1) $K^3\ge \frac{1}{3}$ and this bound is sharp.

(2) If $\frac{1}{3}\le K^3<\frac{5}{14}$, then $p_g(V)=2$,
$n_0(V)=8$, $V$ is canonically fibred by surfaces with $(c_1^2,
p_g)=(1,2)$, and the 7-canonical map of $V$ is generically finite
of degree 2.
\end{thm}

The proof of the main theorem has interesting applications to the
effect that we are able to classify some classes of 3-folds up to
explicit structures. For instance, we have

\begin{thm}\label{bp} Let $V$ be a smooth 3-dimensional projective variety of general
type. The following is true:

(1) if $p_g(V)\ge 3$, then $K^3\ge 1$ and this bound is sharp;

(2) if $p_g(V)\ge 4$, then $K^3\ge 2$ and this bound is sharp;

(3) if $K^3<\frac{1}{2}$ and $p_g(V)=2$, then $V$ is canonically
fibred by surfaces with $c_1^2=1$. In particular, $q(V)=0$ and
$h^2({\mathcal O}_V)\le 1$;

(4) if $K^3< \frac{4}{3}$ and $p_g(V)=3$, then $V$ is canonically
fibred by curves of genus 2 over a birationally ruled surface and
the 4-canonical map is generically finite of degree 2.

(5) if $K^3< \frac{9}{4}$ and $p_g(V)=4$, then $V$ is either a
double cover over $\bbP^3$ or canonically fibred by curves of
genus 2 over a ruled surface.
\end{thm}

All the above statements have supporting examples. More detailed
classification to those $V$ with $p_g=2$ and $3$ is given in
section 2 and section 3. The method works for the case $p_g\le 1$
for which we omit the effective result believing that it might be
far from sharp. Instead we present effective results on 4-folds in
the last section.

\begin{thm}\label{4-folds}
Let $Y$ be a minimal (i.e. $K_Y$ being nef) projective 4-fold of
general type with only canonical singularities. Assume $p_g(Y)\ge
2$ and that $Y$ is not canonically fibred by 3-folds of geometric
genus 1. Then $n_0(Y)\le 24$. Furthermore $K_Y^4\ge \frac{1}{81}$
if $Y$ is not canonically fibred by any irrational pencil of
3-folds either.
\end{thm}
\bigskip

{\bf Notations}

\begin{tabbing}
 \= aaaaaaaaaaaaaaaaa \= bbbbbbbbbbbbbbbbbbbbbbbbbbbbbbb \kill

\> $K^3$    \> the canonical volume of a 3-fold in question\\
\> $p_g$    \> the geometric genus\\
\> $q(V)$    \> the irregularity of $V$\\
\> $h^2({\mathcal O}_V)$ \> the second irregularity of a 3-fold
$V$\\
\> $\chi({\mathcal O})$ \> the Euler Poincare characteristic\\
\> $(c_1^2, p_g)$ \> invariants of a minimal surface of general
type\\
\> $b:=g(B)$ \> the genus of a curve $B$\\
\> $=_{\bQ}$  \> ${\bQ}$-linear equivalence\\
\>$\equiv$     \> numerical equivalence\\
\> $\sim$      \> linear equivalence\\
\> $\Phi_{|L|}$  \> the rational map corresponding to the \\
\>   \   \> linear
system $|L|$\\
\> $\roundup{\cdot}$  \> the round up of a ${\mathbb Q}$-divisor\\
\>$\varphi_m$ \> the m-th pluricanonical map\\
\> $n_0(Y)$ \> the canonical stability index of $Y$\\
\> $D|_S$ \> the restriction of the divisor $D$ to $S$\\
\> $|D||_S$ \> the restriction of the linear system $|D|$ to $S$\\
\> $P_m(V)$  \> the m-th plurigenus of $V$\\
\> $D\cdot C$ \> the intersection number between a \\
\> \  \>divisor D and
curve $C$\\
\> $\pi:X'\rightarrow X$  \> a smooth birational modification\\
\> $f:X'\rightarrow B$  \> an induced fibration from $\varphi_1$\\
\> $M$ \> the movable part of $|K_{X'}|$\\
\> $Z$ \> the fixed part of $|K_{X'}|$\\
\> $S$ \>  a generic irreducible element of $|M|$\\
\> $\sigma:S\rightarrow S_0$ \> contraction onto the minimal
model\\
\> $|G|$ \> a movable linear system on $S$\\
\> $C$ \> a generic irreducible element of $|G|$\\
\> $\xi$ \> the intersection number $\pi^*(K_X)\cdot C$ on $X'$\\
\> $\beta$ \> a positive real number such that\\
\> \   \> $\pi^*(K_X)|_S-\beta C$ is pseudo effective\\
\> $M_i$ \> the movable part of $|iK_{X'}|$\\
\> $S_2$ \> a generic irreducible element of $|M_2|$\\
\> $L_2$ \> the divisor $M_2|_{S_2}$ on $S_2$\\
\end{tabbing}

\section{\bf The case $p_g(V)=2$}
\begin{setup}\label{symbols}{\bf Setting.} Given a smooth
projective 3-dimensional variety of general type with $p_g(V)\ge
2$, let us consider a minimal model as explained in the
introduction. So we assume that $X$ is a minimal model of $V$ with
at worst $\bQ$-factorial terminal singularities. We know
$p_g(X)=p_g(V)$. Take a birational modification $\pi:
X'\rightarrow X$, which exists by Hironaka's theorem, such that

(i) $X'$ is smooth;

(ii) the movable part of $|K_{X'}|$ is base point free.

(iii) $\pi^*(K_X)$ can be written as an effective $\bQ$-divisor
with normal crossings.

Denote by $g$ the composition $\varphi_{1}\circ\pi$. So $g:
X'\longrightarrow W\subseteq{\mathbb P}^{p_g(X)-1}$ is a morphism.
Let $g: X'\overset{f}\longrightarrow B\overset{s}\longrightarrow
W$ be the Stein factorization of $g$. So we have the following
commutative diagram:
$$\begin{CD}
X' @>{f}>> B\\
@V{\pi}VV   @VV{s}V\\
X  @>>\varphi_{1}> W
\end{CD}$$

Write
$$K_{X'}=_{\Bbb Q}\pi^*(K_X)+E=M+Z,$$
where $M$ is the movable part of $|K_{X'}|$, $Z$ the fixed part
and $E$ an effective ${\Bbb Q}$-divisor which is a ${\Bbb Q}$-sum
of distinct exceptional divisors. Throughout we always mean
$\pi^*(K_X)$ by $K_{X'}- E$. So, whenever we take the round up of
$m\pi^*(K_X)$, we always have $\roundup{m\pi^*(K_X)}\le mK_{X'}$
for all positive number $m$. We may also write
$$\pi^*(K_X)=_{\Bbb Q} M+E',$$
where $E'=Z-E$ is actually an effective ${\Bbb Q}$-divisor.

When $\dim\varphi_{1}(X)=2$, we see that a general fiber of $f$ is
a smooth projective curve $C$ of genus $g\ge 2$. When
$\dim\varphi_{1}(X)=1$, we may see from the subadditivity (1.7) in
Mori's paper \cite{Mori} that a general fiber $S$ of $f$ is a
smooth projective surface $S$ of general type. The invariants of
$S$ are $(c_1^2,p_g)=(c_1^2(S_0), p_g(S))$
 where $S_0$ is the minimal model of $S$ and $\sigma: S\longrightarrow S_0$
 the contraction.

We always mean {\it a generic irreducible element $S$ of} $|M|$ by
either a general member of $|M|$ whenever $\dim\varphi_{1}(X)\ge
2$ or, otherwise, a general fiber of $f$.

\begin{defn}
By abuse of concepts, we will also define {\it a generic
irreducible element} $S$ of an arbitrary linear system $|M|$ on a
general variety $V$ in a similar way. Assume that $|M|$ is
movable. A generic irreducible element $S$ is defined to be a
generic irreducible component in a general member of $|M|$. So if
$|M|$ is composed with a pencil (i.e. $\dim\Phi_{|M|}(V)=1$),
$S\le M$ and $M\equiv tS$ for some integer $t\ge 1$. Clearly it
may happen that sometimes $S\not\sim M$.
\end{defn}

We need the following lemma in the proof.

\begin{lem}\label{replace} Keep the same notation as above.
Let $X$ be a minimal 3-fold of general type with at worst
$\bQ$-factorial terminal singularities. Let $\pi: X'\rightarrow X$
be the same modification as above. Suppose we are in the situation
$\dim\varphi_{1}(X)=1$. Let $f:X'\longrightarrow B$ be an induced
fibration. Suppose $g(B)>0$. Then $X$ has another minimal model
$Y$ such that there is a fibration $f_Y:Y\longrightarrow B$ which
is induced by the canonical pencil $|K_Y|$. In particular, the
movable part of $|K_Y|$ is base point free. (Therefore, sometimes,
we will replace $X$ by $Y$ and we simply study on $Y$.)
\end{lem}
\begin{proof}
According to the MMP (see \cite{KMM, K-M, Reid83}), we may take a
relatively minimal model $f_Y: Y\longrightarrow B$ of $f$ such
that $K_Y$ is $f_Y$-nef and that $Y$ is $\bQ$-factorial with
terminal singularities. According to Theorem 1.4 of \cite{Ohno},
$K_{Y/B}$ is nef. (This is a direct consequence of the well-known
semi-positivity of higher direct image of dualizing sheaves by
Fujita, Kawamata, Koll\'ar, Nakayama and Viehweg.) Because
$g(B)>0$, $K_Y$ is nef. Thus $Y$ is a minimal model of $X'$. Take
a common smooth birational modification $X''$ over both $X'$ and
$Y$. We then have the following commutative diagram:
$$\begin{CD}
X'' @>{\theta}>> Y\\
@V{{\theta}'}VV   @VV{f_Y}V\\
X'  @>> f > B
\end{CD}$$

We know that the movable part of $K_{X''}$ is base point free.
Because $p_g(X'')=p_g(Y)$ and $\theta$ is a birational morphism,
we see that the push forward of the movable part of $|K_{X''}|$
onto $Y$ is the movable part of $|K_Y|$. We are done.
\end{proof}
\end{setup}

\begin{setup}\label{estimate}{\bf A known theorem.} In order to
frequently apply it, we rephrase a very effective method (Theorem
2.2 in \cite{IJM}) on how to estimate certain intersection numbers
on $X$ as follows:
\begin{quote}
Let $X$ be a minimal projective 3-fold of general type with only
${\Bbb Q}$-factorial terminal singularities and assume $p_g(X)\ge
2$. Keep the same notations as in \ref{symbols}. Pick up a generic
irreducible element $S$ of $|M|$. Suppose that, on the smooth
surface $S$, there is a movable linear system $|G|$ and set $C$ to
be a generic irreducible element of $|G|$. Denote
$\xi:=(\pi^*(K_X)\cdot C)_{X'}$ and
$$p:=\begin{cases} 1  &\text{if}\ \dim\varphi_1(X)\ge 2\\
p_g(X)-1  &\text{otherwise.}
\end{cases}$$

Assume

\item $\bullet$\  (i) there is a positive integer $m$ such that
the linear system
$$|K_S +\roundup{(m-2)\pi^*(K_X)|_S}|$$
(of which the corresponding rational map is denoted by $\Phi$)
separates different generic irreducible elements of $|G|$. Namely,
if $C_1$ and $C_2$ are different generic irreducible elements of
$|G|$, then $\overline{\Phi(C_1)}\neq \overline{\Phi(C_2)}$;

\item $\bullet$\  (ii) there is a rational number $\beta>0$ such
that $\pi^*(K_X)|_S-\beta C$ is numerically equivalent to an
effective ${\Bbb Q}$-divisor;

\item $\bullet$\  (iii) either the inequality
$$\alpha_0:=(\roundup{(m-1-\frac{1}{p}-\frac{1}{\beta})\pi^*(K_X)|_S}\cdot C)_S\ge 2$$
or at least
 $\alpha:=(m-1-\frac{1}{p}-\frac{1}{\beta})\xi>1$ holds (observe in fact that $\alpha_0\ge \alpha$).
\medskip

 Then we have the inequality $m\xi\ge 2g(C)-2+\alpha_0$.
Furthermore, $\varphi_m$ of $X$ is birational onto its image
provided either $\alpha>2$ or $\alpha_0=2$ and $C$ is
non-hyperelliptic.
\end{quote}
\end{setup}

\begin{rem} As far as the situation is involved, we don't have to verify the
condition (i) in this paper since that was already done in
\cite{IJM}. The main reason is that we are treating the case with
$p_g(X)>0$. On how to take a $\beta$ satisfying the assumption
\ref{estimate}(ii), please refer to \ref{takebeta} below.
\end{rem}

\begin{setup}\label{2K}{\bf Base point freeness of $|2K|$.}
According to Bombieri \cite{Bom}, Reider \cite{Reider},
Catanese-Ciliberto \cite{C-C} and P. Francia \cite{Francia} (or
directly referring to Theorem 3.1 in the survey article by
Ciliberto \cite{Ci}), the bicanonical system $|2K|$ of a minimal
surface of general type with $p_g>0$ is base point free. We will
frequently apply this important result in the context.
\end{setup}

\begin{setup}\label{nt}{\bf Numerical type of a surface of general type.}
Given a smooth projective surface $S$ of general type, we denote
by $\sigma:S\longrightarrow S_0$ the contraction onto the minimal
model. For the need of our proof and according to the standard
surface theory (see \cite{BPV} and \cite{Bom}), we classify $S$
into the five numerical classes as follows:

(1) $S$ is of type $(1,1)$ if $K_{S_0}^2=1$ and $p_g(S)=1$;

(2) $S$ is of type $(1,2)$ if $K_{S_0}^2=1$ and $p_g(S)=2$;

(3) $S$ is of type $(2,3)$ if $K_{S_0}^2=2$ and $p_g(S)=3$;

(4) $S$ is of type $2^+$ if $K_{S_0}^2\ge 2$ but $S$ is not of
type $(2,3)$ (all surfaces with $p_g\ge 3$ fall into this type
because of the Noether inequality: $K^2\ge 2p_g-4$);

(5) $S$ is of type $(1,0)$ if $K_{S_0}^2=1$ and $p_g(S)=0$.

Because $p_g(V)\ge 2$, we see that $p_g(S)>0$. So type $(1,0)$
never happens under the assumption $p_g(V)\ge 2$.
\end{setup}

\begin{setup}\label{takebeta}{\bf How to take a suitable $\beta$?}
We keep the notations in both \ref{symbols} and \ref{estimate}.
Assume $p_g(X)\ge 2$, $\dim\varphi_1(X)=1$ and $b=g(B)=0$. Let $S$
be a general fiber of the induced fibration $f:X'\longrightarrow
B$. Denote by $\sigma:S\longrightarrow S_0$ the contraction onto
the minimal model. By section 4 (at page 526 and page 527) of
\cite{IJM}, we may always choose a sequence of rational numbers
$\beta_0\mapsto\frac{1}{2}$ with $\beta_0<\frac{1}{2}$ such that
$\pi^*(K_X)|_S-\beta_0\sigma^*(K_{S_0})$ is numerically equivalent
to an effective ${\mathbb Q}$-divisor. So we may take a suitable
$\beta$ accordingly:
\begin{quote}
\item $\bullet$ If $p_g(S)>0$ and take $G:=2\sigma^*(K_{S_0})$, we
may choose $\beta=\frac{1}{2}\beta_0$;

\item $\bullet$ If $S$ is of type (1,2) and take $G$ to be the
movable part of $|\sigma^*(K_{S_0})|$, we may choose
$\beta=\beta_0$;

\item $\bullet$ If $S$ is of type (2,3) and take $G$ to be the
movable part of $|\sigma^*(K_{S_0})|$, we may choose
$\beta=\beta_0$.
\end{quote}
\medskip

Usually $\beta$ can be larger whenever $p_g(X)$ is larger. This is
clear by virtue of the method in \cite{IJM}. For example, if
$p_g(X)\ge 3$, we can find a larger $\beta$ as we will do in the
next section.
\end{setup}

\begin{setup}\label{set}{\bf Set up for the case $p_g(V)=2$.}
{}From now on within this section, we assume $p_g(V)=2$. Then $B$
is a smooth curve. Let $S$ be a general fiber of $f$. One may
always write $M=\sum_{i=1}^{a_1}S_i$ where all those $S_i's$ are
disjoint from each other and $a_1\ge p_g(V)-1$. So we have
\begin{align*}
\tag{2.1} K_X^3=(\pi^*(K_X))^3=a_1\pi^*(K_X)^2\cdot
S+\pi^*(K_X)^2\cdot E'.
\end{align*}
Assume $b=g(B)=0$. We have $p=1$ by definition. On the other hand,
we set the divisor $G$ on $S$ as follows:
$$
G:= \begin{cases} \text{the movable part $C$ of
$|\sigma^*(K_{S_0})|$,}
&\text{if $S$ is of type $(1,2)$ or $(2,3)$}\\
2\sigma^*(K_{S_0}), & \text{if $S$ is of type $2^+$ or $(1,1)$.}
\end{cases}$$
Also set
$$\xi:=(\pi^*(K_X)|_S\cdot G)_S=(\pi^*(K_X)\cdot G)_{X'}=
(K_X\cdot \pi_*(G))_X$$ which is independent of the modification
$\pi$ by the intersection theory. All these settings are to
prepare for estimating $\xi$ by means of the technique
\ref{estimate}.
\end{setup}

\begin{setup}\label{2b>0}{\bf The case $b>0$.}
Recall the nefness of $K_X$. This case is simple since the movable
part of $|K_X|$ is already base point free (see Lemma
\ref{replace} for a possible replacement of $X$). So one has
$$\pi^*(K_X)|_S\sim \sigma^*(K_{S_0})$$
and $a_1\ge p_g(V)=2$. Thus we have
$$K_X^3\ge
2\sigma^*(K_{S_0})^2\ge 2.$$ Next we only have to study the case
$b=0$.
\end{setup}

\begin{setup}\label{22^+}{\bf The type $2^+$.}
Because $\pi^*(K_X)$ is nef and big, the equation (2.1) in
\ref{set} gives $K_X^3\ge \pi^*(K_X)^2\cdot S$. So the main point
is to estimate $\pi^*(K_X)^2\cdot S$. By \ref{takebeta}, one may
find a sequence of rational numbers $\beta>0$ with
$\beta\rightarrow \frac{1}{4}$ and $\beta<\frac{1}{4}$ such that
$\pi^*(K_X)|_S-\beta G$ is numerically equivalent to an effective
$\bQ$-divisor. Thus we have
$$\pi^*(K_X)^2\cdot S\ge \beta \xi.$$
We are reduced to estimate the rational number $\xi$. Let $C$ be a
generic irreducible element of $|G|$.  We know that $|G|$ is base
point free. So $C$ is a smooth curve with $\deg(K_C)\ge 12$. Note
that $C$ must be non-hyperelliptic since $n_0(S)\le 3$ by Bombieri
\cite{Bom} and Reider \cite{Reider}.

Take $m_1=7$. Then
$\alpha=(m_1-1-\frac{1}{p}-\frac{1}{\beta})\xi>0$. So
$\alpha_0>0$. Because in this case $G$ is an even divisor by
definition, we must have $\alpha_0\ge 2$. So \ref{estimate} gives
$\xi\ge 2$. Thus we have $K_X^3\ge \frac{1}{2}$.
\end{setup}

\begin{setup}\label{2(2,3)}{\bf The type $(2,3)$.}
Let $C$ be a generic irreducible element of $|G|$. Then $C$ is a
curve of genus 3 and $C^2=K_{S_0}^2=2$ (see pages 226-227 of
\cite{BPV}). According to \ref{takebeta}, we may take a rational
number $\beta\mapsto \frac{1}{2}$ such that $\pi^*(K_X)|_S-\beta
C$ is pseudo-effective. Thus $\xi\ge\beta C^2\ge 1$. So we have
$$K_X^3\ge \pi^*(K_X)^2\cdot S\ge \beta\xi\ge \frac{1}{2}.$$
\end{setup}

\begin{setup}\label{2(1,1)}{\bf The type $(1,1)$.}
This situation has to be studied in an alternative way. We have an
induced fibration $f:X'\longrightarrow B$. Because $|K_{X'}|$ is
composed with a pencil, $b=0$ and $f_*\omega_{X'}$ is a line
bundle in this case, we see that $\deg f_*\omega_{X'}>0$ by the
Riemann-Roch theorem. According to Bombieri \cite{Bom}, a surface
$S$ of type (1,1) has $q(S)=0$. So $R^1f_*\omega_{X'}=0$. Applying
Koll\'ar's formulae (Proposition 7.6 at page 36 of \cite{Kol}),
one gets $h^2(\so_{X})=h^1(f_*\omega_{X'})=0$ and $q(X)=0$. Thus
$\chi(\so_X)=-1.$ By the plurigenus formula of Reid (see Chapter
III of \cite{YPG}) and omitting the correction term, we have
$$P_2(X)\ge \frac{1}{2}K_X^3-3\chi(\so_X)>3.$$
This means $P_2(X)\ge 4$ because $K_X^3>0$. Let $M_2$ be the
movable part of $|2K_{X'}|$. We consider the natural restriction
map $\gamma$:
$$H^0(X', M_2)\overset{\gamma}\longrightarrow U_2\subset H^0(S, M_2|_S)\subset
H^0(S, 2K_S),$$ where $U_2$ is the image of $\gamma$ as a
$\bC$-subspace of $H^0(S, M_2|_S)$. Because
$h^0(2K_S)=K^2+\chi({\mathcal O}_S)=3$, we see that
$1\le\dim_{\mathbb C}U_2\le 3$. Denote by $\Lambda_2$ the linear
system corresponding to $U_2$. We have
$\dim\Lambda_2=\dim_{\bC}U_2-1$.
\medskip

Case 1. $\dim_{\bC}U_2=3$.

Since $\Lambda_2$ is a sub-system of $|2K_S|$, we see that the
restriction of $\phi_{2,X'}$ to $S$ is exactly the bicanonical map
of $S$. Because $\phi_{2, S}$ is a generically finite morphism of
degree $4$, $\phi_{2,X'}$ is also a generically finite map of
degree $4$. (Here let me slightly explain the fact about $\phi_{2,
S}$. As I have explained in \ref{22^+}, any minimal surface of
general type with $p_g>0$ has a base point free bicanonical
system. It follows directly that the movable part of $|2K_S|$
defines a generically finite morphism of degree 4. Suitable
references might be \cite{BPV} and \cite{Bom}. )

Let $S_2\in|M_2|$ be a general member. We can further remodify
$\pi$ such that $|M_2|$ is base point free too. Then $S_2$ is a
smooth projective irreducible surface of general type. On the
surface $S_2$, denote $L_2:=S_2|_{S_2}$. $L_2$ is a nef and big
divisor. We have
$$2\pi^*(K_X)|_{S_2}\ge S_2|_{S_2}=L_2.$$
We consider the natural map
$$H^0(X', S_2)\overset{\gamma'}\longrightarrow\overline{U_2}\subset H^0(S_2, L_2),$$
where $\overline{U_2}$ is the image of $\gamma'$. Denote by
$\overline{\Lambda_2}$ the linear system corresponding to
$\overline{U_2}$. Because $\varphi_2$ is a generically finite map
of degree $4$, we see that $|L_2|$ has a sub-system
$\overline{\Lambda_2}$ which gives a generically finite map of
degree $4$. By direct calculation on the surface $S_2$, we have
$$L_2^2\ge 4(\dim_{\Bbb C}\overline{\Lambda_2}-1)\ge 4(P_2(X)-3)\ge
4.$$
Therefore we have
$$K_X^3\ge\frac{1}{8}L_2^2\ge\frac{1}{2}.$$
\medskip

Case 2. $\dim_{\bC}U_2=2$.

In this case, $\dim\varphi_2(S)=1$ and $\dim\varphi_2(X)=2$. We
may further remodify $\pi$ such that $|M_2|$ is base point free.
Taking the Stein factorization of $\Phi_{2K_{X'}}$, we get a
induced fibration $f_2:X'\longrightarrow B_2$ where $B_2$ is a
surface. Let $C$ be a general fiber of $f_2$. We see that $S$ is
naturally fibred by curves with the same numerical type as $C$. On
the surface $S$, we have a free pencil $\Lambda_2\subset |2K_S|$
(here we mean the movable part of $\Lambda_2$ has no base points).
Let $|C_0|$ be the movable part of $\Lambda_2$. Then $h^0(S,
C_0)=2$ because $\dim_{\bC}U_2=2$ in this case. Because $q(S)=0$
(\cite{BPV, Bom}), we see that $|C_0|$ is a rational pencil. So a
general member of $|C_0|$ is an irreducible curve with the same
numerical type with $C$.

\begin{claim}\label{11}
Let $S$ be a surface of type (1,1). Suppose that there is an
effective irreducible curve $C$ on $S$ such that $C\le
\sigma^*(2K_{S_0})$ and $h^0(S, C)=2$. Then $C\cdot
\sigma^*(K_{S_0})\ge 2$.
\end{claim}
\begin{proof} We may assume that $|C|$ is base point free. Otherwise,
we blow-up $S$ at base points of $|C|$ and consider the movable
part. Denote $C_1:=\sigma(C)$. Then $h^0(S_0, C_1)\ge 2$. Suppose
$C\cdot \sigma^*(K_{S_0})=1$. Then $C_1\cdot K_{S_0}=1$. Because
$p_a(C_1)\ge 2$, we see that $C_1^2>0$. {}From
$K_{S_0}(K_{S_0}-C_1)=0$, we get $(K_{S_0}-C_1)^2\le 0$, i.e.
$C_1^2\le 1$. Thus $C_1^2=1$ and $K_{S_0}\equiv C_1$. This means
$K_{S_0}\sim C_1$ by virtue of \cite{BPV, Bom, Catan}, which is
impossible because $p_g(S)=1$.  So $C\cdot\sigma^*(K_{S_0})\ge 2$
\end{proof}

According to Claim \ref{11}, we have
$(C_0\cdot\sigma^*(K_{S_0}))_S\ge 2$. We also have
$$(\pi^*(K_X)\cdot C)_{X'}=(\pi^*(K_X)|_S\cdot C_0)_S.$$

According to \ref{takebeta}, there is a rational number
$\beta_0\mapsto \frac{1}{2}$ such that
$\pi^*(K_X)|_S-\beta_0\sigma^*(K_{S_0})$ is pseudo-effective. Thus
$(\pi^*(K_X)|_S\cdot C_0)_S\ge \frac{1}{2}\sigma^*(K_{S_0})\cdot
C_0\ge 1$. Now we study on the surface $S_2$. We may write
$$S_2|_{S_2}\sim\sum_{i=1}^{a_2}C_i\equiv a_2C,$$
where the $C_i's$ are fibers of $f_2$ and $a_2\ge P_2(X)-2\ge 2$.
Noting that $2\pi^*(K_X)|_{S_2}\ge S_2|_{S_2}$, we get
\begin{align*}
4K_X^3&\ge 2\pi^*(K_X)^2\cdot S_2=2(\pi^*(K_X)|_{S_2})_{S_2}^2\\
 &\ge 2(\pi^*(K_X)|_{S_2}\cdot C)_{S_2}\\
 &=2(\pi^*(K_X)\cdot C)_{X'}\ge 2.
\end{align*}
So we still have $K_X^3\ge\frac{1}{2}.$
\medskip

Case 3. $\dim_{\bC}U_2=1$.

In this case,  $\dim\varphi_2(X)=1$. Because $p_g(X)>0$, we see
that both $\varphi_2$ and $\varphi_1$ induce the same fibration
$f:X'\longrightarrow B$ after taking the Stein factorization of
them. So we may write
$$
2\pi^*(K_X)\sim \sum_{i=1}^{a_2'}S_i+E_2' \equiv a_2'S+E_2',$$
where the $S_i's$ are fibers of $f$, $E_2'$ is an effective
$\bQ$-divisor, $a_2'\ge P_2(X)-1\ge 3$ and $S$ is a surface of
type (1,1). So we get
$$2K_X^3\ge 3(\pi^*(K_X)|_S)_S^2\ge 3\beta_0 \pi^*(K_X)|_S\cdot\sigma^*(K_{S_0})$$
where $\beta_0\mapsto \frac{1}{2}.$

Now we may apply \ref{estimate} to estimate
$\pi^*(K_X)|_S\cdot\sigma^*(K_{S_0})$. We have $p=1$ and set
$G:=2\sigma^*(K_{S_0})$. \ref{takebeta} allows us to take
$\beta=\frac{1}{2}\beta_0$. Note that a generic irreducible
element $C\in |G|$ on $S$ is a non-hyperelliptic curve (since we
know $n_0(S)\le 3$ by Bombieri \cite{Bom}) and $\deg(K_C)=6$.
Because $C$ is actually an even divisor on the surface $S$, the
number $\alpha_0$ in \ref{estimate} is a positive even number
whenever $l=m-1-\frac{1}{p}-\frac{1}{\beta}>0$ for some integer
$m$. This means we have $\alpha_0\ge 2$ whenever $l>0$. Now taking
$m_1=7$, we may verify $l>0$ and so all assumptions in
\ref{estimate} are satisfied. Then \ref{estimate} gives
$\xi:=\pi^*(K_X)|_S\cdot G\ge \frac{8}{7}$. Take $m_2=8$. Then
similarly we have $\xi\ge \frac{5}{4}$. We may check that this is
the best bound we could get from this technique. Thus we have
$\pi^*(K_X)|_S\cdot\sigma^*(K_{S_0})\ge \frac{5}{8}$ and so
$K_X^3\ge \frac{15}{32}.$
\medskip

Thus we have $K_X^3\ge \frac{15}{32}$.
\end{setup}

\begin{setup}\label{2(1,2)}{\bf The type $(1,2)$.}
According to \ref{takebeta}, we may take $\beta\mapsto
\frac{1}{2}$. {}From section 4 (page 527) of \cite{IJM}, we know
$\xi\ge \frac{3}{5}$. Take $m_1=8$. Then
$\alpha=(m_1-1-1-\frac{1}{\beta})\xi>2$. So \ref{estimate} gives
$\xi\ge \frac{5}{8}$. Take $m_2=9$. Then
$\alpha=(m_2-2-\frac{1}{\beta})\xi>3$. \ref{estimate} gives
$\xi\ge \frac{2}{3}$. One may check that this is the best we could
get. Thus we have
$$K_X^3\ge \pi^*(K_X)^2\cdot S\ge \beta \xi\ge \frac{1}{3}. $$

We go on studying this case from the point of view of the
canonical stability index.
\end{setup}

\begin{claim}\label{iff}
Under the situation \ref{2(1,2)}, $n_0(V)=8$ if and only if
$\xi=\frac{2}{3}$. In particular $n_0(V)=8$ whenever
$K^3<\frac{5}{14}$.
\end{claim}
\begin{proof}
If $\xi>\frac{2}{3}$, let us take $m=7$. Then
$\alpha=(m-2-\frac{1}{\beta})\xi>2$. \ref{estimate} says that
$\varphi_7$ is birational and $\xi\ge \frac{5}{7}$. So we get
$K_X^3\ge \frac{5}{14}$. This means that $n_0(V)=8$ implies
$\xi=\frac{2}{3}$.

To the contrary, assume $\xi=\frac{2}{3}$, we want to show that
$n_0(V)=8$ and that the 7-canonical map is generically finite of
degree 2.

We consider the sub-system
$$|K_{X'}+\roundup{5\pi^*(K_X)}+S|\subset |7K_{X'}|.$$
This system obviously separates different generic irreducible
elements of $|M|$ because $K_{X'}+\roundup{5\pi^*(K_X)}+S\ge S$
and $|S|$ is composed with a rational pencil. By the Tankeev
principle for birationality, it is sufficient to prove that
$|7K_{X'}||_S$ gives a birational map. Noting that $5\pi^*(K_X)$
is nef and big, the Kawamata-Viehweg vanishing theorem ( see
\cite{E-V, Ka, V}) gives the surjective map
\begin{eqnarray*}
&H^0(X',K_{X'}+\roundup{5\pi^*(K_X)}+S)\\
\longrightarrow  &H^0(S, K_S+\roundup{5\pi^*(K_X)}|_S).
\end{eqnarray*}
We are reduced to prove that $|K_S+\roundup{5\pi^*(K_X)}|_S|$
gives a birational map. We still apply the Tankeev principle.
Because
$$ K_S+\roundup{5\pi^*(K_X)}|_S\ge
K_S+\roundup{5\pi^*(K_X)|_S},$$ the linear system
$|K_S+\roundup{5\pi^*(K_X)}|_S|$ separates different irreducible
elements of $|G|$ where $G$ is defined as in \ref{set}. This is
because $\roundup{5\pi^*(K_X)}|_S\ge 0$ and $K_S\ge C$ and $|C|$
is composed with a rational pencil of curves. Note also here that
$q(S)=0$ and $C$ is a smooth curve of genus 2 (see page 225 of
\cite{BPV}). Now pick up a generic irreducible element $C\in |G|$.
By \ref{takebeta} or \cite{IJM}, there is a rational number
$\beta\mapsto\frac{1}{2}$ and an effective ${\Bbb Q}$-divisor $H$
on $S$ such that
$$\frac{1}{\beta}\pi^*(K_X)|_S\equiv C+H.$$
By the vanishing theorem, we have the surjective map
$$H^0(S, K_S+\roundup{5\pi^*(K_X)|_S-H})\longrightarrow
H^0(C, D),$$ where $D:=\roundup{5\pi^*(K_X)|_S-C-H}|_C$ is a
divisor on $C$. Noting that
$$5\pi^*(K_X)|_S-C-H\equiv (5-\frac{1}{\beta})\pi^*(K_X)|_S$$
and that $C$ is nef on $S$, we have $\deg(D)\geq\alpha$ and thus
$\deg(D)\geq \alpha_0$.

Now if $\xi=\frac{2}{3}$, then $\deg(D)\ge 2$. This means
$|K_C+D|$ is base point free. Noting that $C$ is a curve of genus
2, $|K_C+D|$ gives a finite map of degree $\le $2. Thus
$\varphi_7$ must be a generically finite map if it is not
birational.

Now we show that the 7-canonical map is not birational. Denote by
$M_7$ the movable part of $|7K_{X'}|$ and by $M_7'$ the movable
part of $|K_{X'}+\roundup{5\pi^*(K_X)}+S|$. Apparently, one has
$$M_7'|_S\le M_7|_S\le 7\pi^*(K_X)|_S.$$
Because $7\pi^*(K_X)|_S\cdot C=\frac{14}{3}<5$, we have
$M_7|_S\cdot C\le 4$. On the other hand, we have $M_7'|_S\cdot
C\ge \deg(K_C+D)\ge 4$. This means that $\deg(K_C+D)=4$. It is
obvious that
$$\roundup{7\pi^*(K_X)|_S}|_C\sim 5P$$
where $P$ is a point of $C$ such that $\so_C(2P)\cong \omega_C$.
Furthermore one sees that $(\rounddown{7\pi^*(K_X)}|_S)|_C\sim nP$
with $n\le 4$. Because $M_7\le \rounddown{7\pi^*(K_X)}$, we must
have $n=4$ and $K_C+D\sim 4P$. So $|K_C+D|$ must give a finite map
of degree 2. So is $\varphi_7$.
\end{proof}

\begin{cor}\label{c21} Assume $p_g(V)=2$. If $K^3<\frac{5}{14}$,
then $V$ is of type $(1,2)$ and $n_0(V)=8$.
\end{cor}
\begin{proof}
This is clear by virtue of \ref{2b>0}, \ref{22^+}, \ref{2(2,3)},
\ref{2(1,1)}, \ref{2(1,2)} and \ref{iff}.
\end{proof}

\begin{cor}\label{c22} Assume $p_g(V)=2$. If $K^3<\frac{1}{2}$,
then $V$ is canonically fibred by surfaces of type $(1,1)$ or
$(1,2)$. In particular, $q(V)=0$ and $h^2(\so_V)\le 1$.
\end{cor}
\begin{proof}
This is clear by virtue of \ref{2b>0}, \ref{22^+}, \ref{2(2,3)},
\ref{2(1,1)} and \ref{2(1,2)}.

Note that $q(S)=0$ if $S$ is either of type $(1,2)$ or of type
$(1,1)$ (by Bombieri \cite{Bom}). So $q(V)=0$ by Koll\'ar's
formula.

If $S$ is of type $(1,1)$, then $f_*\omega_{X'}$ is a line bundle
of positive degree because $|K_{X'}|$ is composed with a pencil
and $p_g(X)>1$. So one has, by Koll\'ar's formula (Proposition 7.6
of \cite{Kol}),
$$h^2(\so_V)=h^2(\so_X)=h^1(B,f_*\omega_{X'})+h^0(B,
R^1f_*\omega_{X'})=0.$$

If $S$ is of type $(1,2)$, according to the result (the middle
part) at page 524 of \cite{IJM}, we have $h^2(\so_V)\le 1$. We are
done.
\end{proof}

There is an example which shows that our estimation is sharp.
\begin{exmp}\label{1/3} In \cite{C-R}, Fletcher found a canonical
3-fold with $K^3=\frac{1}{3}$ and $p_g=2$ as a hypersurface
$X_{16}\subset \bbP(1,1,2,3,8)$. The example has 3 terminal
singularities of type $2\times \frac{1}{2}(1,-1,1)$, $1\times
\frac{1}{3}(1,-1,1)$. This fits in with our argument above since
this 3-fold has the canonical stability index 8, and the
7-canonical map is generically finite of degree 2.
\end{exmp}

\begin{rem} Probably there are more examples.
It would be very interesting to find new examples with $p_g=2$ and
$n_0=8$.
\end{rem}

\section{\bf The case $p_g(V)\ge 3$ and the proof of Theorem \ref{main}}
Assume from now on that $p_g(V)\ge 3$. Set $d:=\dim \varphi_1(X)$.

\begin{setup}\label{d3}{\bf The case $d=3$.}
In this case, $p_g(X)\ge 4$. It is obvious that $K^3\ge 2$. In
fact, one has a general inequality $K^3\ge 2p_g(X)-6$ which may be
obtained by induction on the dimension ( see the main theorem in
\cite{Kob}). Anyway that is another kind of question.
\end{setup}

\begin{setup}\label{d2}{\bf The case $d=2$.}
We have an induced fibration $f:X'\longrightarrow B$ from
$\varphi_1$. By an easy subadditivity (1.7 in \cite{Mori}), a
general fibre $C$ is a smooth curve of genus $\ge 2$. We have
$\pi^*(K_X)=S+E'$ and then
$$\pi^*(K_X)|_S=S|_S+E'|_S.$$
Because $d=2$, $S|_S\equiv a_2C$ where $a_2\ge p_g(X)-2$ and that
the equality holds if and only $|S|_S|$ is a rational pencil. Thus
$$K^3\ge \pi^*(K_X)^2\cdot S\ge
a_2\pi^*(K_X)\cdot C.$$

We have $p=1$ and may set $\beta=a_2\ge 1$. Set $G:=C$. We hope to
run \ref{estimate}. In the proof of Case 2 of Theorem 3.1 at page
521 in \cite{IJM}, we have shown $\xi\ge \frac{6}{7}$. Take
$m_1=8$. Then $\alpha=(m_1-1-1-\frac{1}{\beta})\xi>4$. This gives
$\xi\ge \frac{7}{8}$. Similarly one may use induction to see
$\xi\ge \frac{k}{k+1}$ for all $k>8$. Taking limits, one gets
$\xi\ge 1$.

So when $p_g(X)=3$, we have $K_X^3\ge 1$.

When $p_g(X)\ge 4$, we have $K_X^3\ge 2$.

In fact, if $g(C)\ge 3$ and take $m=6$, the technique
\ref{estimate} gives $\xi>1$. Whenever the equalities for $K^3$
hold above, one must have $\xi=1$. So we see that $C$ must be a
curve of genus 2, and that $|S|_S|$ must be a rational pencil (by
Riemann-Roch) and so that $B$ is a birationally ruled surface.
\end{setup}

\begin{setup}\label{d1>}{\bf The case $d=1$ and $b>0$.}
In this case, one has $\pi^*(K_X)\equiv \bar{a_2}S+E'$ where
$\bar{a_2}\ge p_g(X)\ge 3$. Also one has $\pi^*(K_X)\sim
\sigma^*(K_{S_0})$ because one may take a $X$ such that the
movable part of $|K_X|$ is base point free by Lemma \ref{replace}.
One apparently has $K_X^3\ge 3$. We are done.
\end{setup}

Now we see how to find a bigger $\beta$ than those found in the
last section.

\begin{lem}\label{beta} Assume $p_g(X)=p+1\ge 3$, $d=1$ and $b=0$.
Then one may take a rational number $\beta_0\mapsto\frac{p}{p+1}$
such that $\pi^*(K_X)-\beta_0\sigma^*(K_{S_0})$ is
pseudo-effective.
\end{lem}
\begin{proof}
Because $|K_X|$ is composed with a rational pencil and $P_g(X)=
p+1$, one has
$$\so_B(p)\hookrightarrow {f}_*\omega_{X'}.$$
Thus we have
$${f}_*\omega_{X'/B}^{4p}\hookrightarrow
{f}_*\omega_{X'}^{4p+8}.$$

For any integer $k$, denote by $M_k$ the movable part of
$|kK_{X'}|$. Note that ${f}_*\omega_{X'/B}^{4p}$ is generated by
global sections (see \cite{F, Ka1, Kol, N, V2}) and so that any
local section can be extended to a global one. On the other hand,
$|4p\sigma^*(K_{S_0})|$ is base point free and is exactly the
movable part of $|4pK_S|$ by Bombieri \cite{Bom} or Reider
\cite{Reider}. Applying Lemma 2.7 of \cite{MPC}, one has
$$(4p+8)\pi^*(K_X)|_S\ge M_{4p+8}|_S\ge 4p\sigma^*(K_{S_0}).$$
This means that there is an effective $\bQ$-divisor $E_0'$ such
that
$$(4p+8)\pi^*(K_X)|_S=_{\bQ} 4p\sigma^*(K_{S_0})+E_0'.$$
Thus
$$\pi^*(K_X)|_S\equiv \frac{p}{p+2}\sigma^*(K_{S_0})+E_0$$
where $E_0=\frac{1}{4p+8}E_0'$ is an effective $\bQ$-divisor. Set
$a_0:=4p+8$ and $b_0:=4p$.

Assume we have defined $a_n$ and $b_n$. We describe $a_{n+1}$ and
$b_{n+1}$ inductively such that $\beta_0\ge
\frac{b_{n+1}}{a_{n+1}}$. One may assume from the beginning that
$a_n\pi^*(K_X)$ supports on a divisor with normal crossings. Then
the Kawamata-Viehweg vanishing theorem implies the surjective map
$$H^0(K_{X'}+\roundup{a_n\pi^*(K_X)}+S)\longrightarrow H^0(S, K_S+
\roundup{a_n\pi^*(K_X)}|_S).$$ That means
\begin{eqnarray*}
|K_{X'}+\roundup{a_n\pi^*(K_X)}+S||_S&=&|K_S+\roundup{a_n\pi^*(K_X)}|_S|\\
&\supset& |K_S+b_n\sigma^*(K_{S_0})|\\
&\supset& |(b_n+1)\sigma^*(K_{S_0})|.
\end{eqnarray*}
Denote by $M_{a_n+1}'$ the movable part of $|(a_n+1)K_{X'}+S|$.
Applying Lemma 2.7 of \cite{MPC} again, one has
$$M_{a_n+1}'|_S\ge (b_n+1)\sigma^*(K_{S_0}).$$
Re-modifying our original $\pi$ such that $|M_{a_n+1}'|$ is base
point free. In particular, $M_{a_n+1}'$ is nef. According to
\ref{kr}, $|mK_X|$ gives a birational map whenever $m\ge 6$. Thus
$M_{a_n+1}'$ is big.

Now the Kawamata-Viehweg vanishing theorem gives
\begin{eqnarray*}
|K_{X'}+M_{a_n+1}'+S||_S&=&|K_S+M_{a_n+1}'|_S|\\
&\supset& |K_S+(b_n+1)\sigma^*(K_{S_0})|\\
&\supset& |(b_n+2)\sigma^*(K_{S_0})|.
\end{eqnarray*}

We may repeat the above procedure inductively. Denote by
$M_{a_n+t}'$ the movable part of $|K_{X'}+M_{a_n+t-1}'+S|$. For
the same reason, we may assume $|M_{a_n+t}'|$ to be base point
free. Inductively one has
$$M_{a_n+t}'|_S\ge (b_n+t)\sigma^*(K_{S_0}).$$
Applying the vanishing theorem once more, we have
\begin{eqnarray*}
|K_{X'}+M_{a_n+t}'+S||_S&=&|K_S+M_{a_n+t}'|_S|\\
&\supset& |K_S+(b_n+t)\sigma^*(K_{S_0})|\\
&\supset& |(b_n+t+1)\sigma^*(K_{S_0})|.
\end{eqnarray*}

Take $t=p$. Noting that
$$|K_{X'}+M_{a_n+p}'+S|\subset |(a_n+p+1)K_{X'}|$$
and applying Lemma 2.7 of \cite{MPC} again, one has
$$(a_n+p+1)\pi^*(K_X)|_S\ge M_{a_n+p+1}|_S\ge
(b_n+p)\sigma^*(K_{S_0}).$$ Set $a_{n+1}:=a_n+p+1$ and
$b_{n+1}=b_n+p$. We have seen
$$\beta_0\ge \frac{b_{n+1}}{a_{n+1}}.$$
 A direct calculation gives
 $$a_n=n(p+1)+a_0$$
$$b_n=np+b_0.$$
Take the limit with $n\mapsto +\infty$, one has $\beta_0\ge
\frac{p}{p+1}$. We are done.
\end{proof}

\begin{setup}\label{32^+}{\bf The type $2^+$.}
We have $K_X^3\ge (p_g(X)-1)\pi^*(K_X)^2\cdot S$ by the equation
(2.1). The main point is to estimate $\pi^*(K_X)^2\cdot S$. We
still set $G=2\sigma^*(K_{S_0})$. We have $p\ge 2$ by definition.
By \ref{beta}, one may find a rational number
$\beta=\frac{1}{2}\beta_0\mapsto\frac{1}{3}$ whenever $p_g(X)=3$
(or $\frac{3}{8}$ whenever $p_g(X)\ge 4$) such that
$\pi^*(K_X)|_S-\beta G$ is pseudo-effective. Now we have
$$\xi=\pi^*(K_X)|_S\cdot G\ge \beta G^2=4K_{S_0}^2\beta\ge 8\beta.$$
Thus we have $K^3\ge (p_g(X)-1)\beta\xi.$

This gives $K^3\ge \frac{16}{9}$ if $p_g(X)=3$ and $K^3\ge
\frac{27}{8}$ if $p_g(X)\ge 4$.
\end{setup}

\begin{setup}\label{3(2,3)}{\bf The type $(2,3)$.}
Let $C$ be a generic irreducible element of $|G|$. Then $C$ is a
curve of genus 3 as we know. It has been shown in \ref{2(2,3)}
that $\xi\ge 1$. On the other hand, \ref{beta} tells that one may
take a rational number $\beta\mapsto\frac{2}{3}$ whenever
$p_g(X)=3$ (or $\frac{3}{4}$ whenever $p_g(X)\ge 4$).

Consider the case $p_g(X)=3$. Take $m_1=5$ and we have $p=2$. Then
$\alpha=(m_1-1-\frac{1}{2}-\frac{1}{\beta})\xi\ge 2$. Then one has
$\xi\ge \frac{6}{5}$. Take $m_2=5$ again. One sees that
$\alpha>2$. This means that $\varphi_5$ is birational. Take
$m_3=4$. Then $\alpha=(4-1-\frac{1}{2}-\frac{1}{\beta})\xi>1$. One
has $\xi\ge \frac{3}{2}$. It seems that this is the best one can
get under the condition $p_g(X)=3$. Thus $K^3\ge 2\beta\xi\ge 2$.

Whenever $p_g(X)\ge 4$, one has $K^3\ge 3$. There might be better
bound.
\end{setup}

\begin{setup}\label{3(1,2)}{\bf The type $(1,2)$.}
This is a continuation of \ref{2(1,2)} under the assumption
$p_g(X)\ge 3$. We have got $\xi\ge \frac{2}{3}$ in \ref{2(1,2)}
because the proof for $p_g(X)=2$ is automatically true for
$p_g(X)\ge 3$. Recall that $G$ is the movable part of
$|\sigma^*(K_{S_0})|$. We have $p=2$ and $\beta=\beta_0$ is near
$\frac{2}{3}$ by \ref{beta}. Take $m=5$. Then
$\alpha=(5-1-\frac{1}{2}-\frac{1}{\beta})\xi>1$. This gives
$\xi\ge \frac{4}{5}$. Take $m=6$. Then it gives $\xi\ge
\frac{5}{6}$. An induction step gives $\xi\ge 1$.

Now if $p_g(X)=3$, then $K^3\ge 2\beta\xi\ge \frac{4}{3}$.

If $p_g(X)=4$, then $K_X^3\ge 3\beta\xi\ge \frac{9}{4}$.
\end{setup}

\begin{setup}\label{3(1,1)}{\bf The type $(1,1)$.}
Comparing with \ref{2(1,1)}, we have better situation because
$p_g(X)\ge 3$. We keep the same notations and pace as in
\ref{2(1,1)}.

We still consider the natural restriction map $\gamma$:
$$H^0(X', M_2)\overset{\gamma}\longrightarrow U_2\subset H^0(S, M_2|_S)\subset
H^0(S, 2K_S),$$ where $U_2$ is the image of $\gamma$ as a
$\bC$-subspace of $H^0(S, M_2|_S)$. If $p_g(X)=3$, then a similar
calculation to that in \ref{2(1,1)} gives $\chi(\so_X)=-2$ and
then the plurigenus formula of Reid gives $P_2(X)\ge 7$. If
$p_g(X)\ge 4$, then $\chi(\so_X)\le -3$ and $P_2(X)\ge 10$.
\medskip

Case 1. $\dim_{\bC}U_2=3$.

All those arguments may be copied from Case 1 of \ref{2(1,1)}. We
only write down the estimation. We have
$$2\pi^*(K_X)|_{S_2}\ge S_2|_{S_2}=L_2.$$
By direct calculation on the surface $S_2$, we have
$$L_2^2\ge 4(\dim_{\Bbb C}\overline{\Lambda_2}-1)\ge 4(P_2(X)-3).$$
Therefore we have
$$K_X^3\ge\frac{1}{8}L_2^2\ge \frac{1}{2}(P_2(X)-3).$$

If $p_g(X)=3$, then $K^3\ge 2$.

If $p_g(X)\ge 4$, then $K^3\ge \frac{7}{2}$.
\medskip

Case 2. $\dim_{\bC}U_2=2$.

We may copy the proof in Case 2 of \ref{2(1,1)}. In this case,
$\varphi_2$ induces a fibration $f_2: X'\longrightarrow B_2$ with
general fibre a smooth curve $C$. Lemma \ref{beta} allows us to
choose a rational number $\beta_0\mapsto\frac{2}{3}$ when
$p_g(X)=3$ (or $\frac{3}{4}$ when $p_g(X)\ge 4$) such that
$\pi^*(K)|_S-\beta_0\sigma^*(K_{S_0})$ is pseudo-effective. Thus
$$(\pi^*(K_X)\cdot C)_{X'}=\pi^*(K_X)|_S\cdot C\ge
\beta_0\sigma^*(K_{S_0})\cdot C\ge 2\beta_0$$ where Claim \ref{11}
is applied once more to get the inequality. Recall that $S_2$ is a
general member of the movable part of $|2K_{X'}|$. We may write
$$S_2|_{S_2}\sim\sum_{i=1}^{a_2}C_i\equiv a_2C,$$
where the $C_i's$ are fibers of $f_2$ and $a_2\ge P_2(X)-2$.
Noting that $2\pi^*(K_X)|_{S_2}\ge S_2|_{S_2}$, we get
\begin{align*}
4K_X^3&\ge 2\pi^*(K_X)^2\cdot S_2=2(\pi^*(K_X)|_{S_2})_{S_2}^2\\
 &\ge (P_2-2)(\pi^*(K_X)|_{S_2}\cdot C)_{S_2}\\
 &=(P_2-2)(\pi^*(K_X)\cdot C)_{X'}\ge 2(P_2-2)\beta_0.
\end{align*}

If $p_g(X)=3$, then $K^3\ge \frac{5}{3}$.

If $p_g(X)\ge 4$, then $K^3\ge 3$.
\medskip

Case 3. $\dim_{\bC}U_2=1$.

Similar to the proof in Case 3 of \ref{2(1,1)}, we have
$$2\pi^*(K_X)\equiv a'S+E'$$
where $a'\ge P_2(X)-1$. So we get
$$2K_X^3\ge (P_2-1)(\pi^*(K_X)|_S)_S^2\ge (P_2-1)\beta_0
\pi^*(K_X)|_S\cdot\sigma^*(K_{S_0})$$ where $\beta_0\mapsto
\frac{2}{3}$ when $p_g(X)=3$  (or $\frac{3}{4}$ otherwise) by
Lemma \ref{beta}.

Now we have to estimate $\pi^*(K_X)|_S\cdot\sigma^*(K_{S_0})$. We
have $p=2$ and set $G:=2\sigma^*(K_{S_0})$ and
$\beta=\frac{1}{2}\beta_0$. Let $C$ be a generic irreducible
element of $|G|$. Then $C$ is non-hyperelliptic as we have seen
and $\deg(K_C)=6$. Take $m_1=6$. Then
$\alpha=(m_1-1-\frac{1}{2}-\frac{1}{\beta})\xi>0$. The even
divisor $G$ gives $\alpha_0\ge 2$. \ref{estimate} gives $\xi\ge
\frac{4}{3}$. This might be the best that we get.

We have $K^3\ge \frac{1}{2}(P_2-1)\beta\xi.$

If $p_g(X)=3$, then $K^3\ge \frac{4}{3}$.

If $p_g(X)\ge 4$, then $K^3\ge \frac{9}{4}$. We are done.

\end{setup}

\begin{cor}\label{sum}
Let $V$ be a smooth projective 3-fold of general type. Then

(1) $K^3\ge 1$ whenever $p_g(V)=3$;

(2) $K^3\ge 2$ whenever $p_g(V)=4$.
\end{cor}
\begin{proof}
This is clear by \ref{d3}, \ref{d2}, \ref{d1>}, \ref{32^+},
\ref{3(2,3)}, \ref{3(1,2)} and \ref{3(1,1)}.
\end{proof}

\begin{setup}{\bf Proof of Theorem \ref{main}.}
\begin{proof}
\ref{2b>0}, \ref{22^+}, \ref{2(2,3)}, \ref{2(1,1)}, \ref{2(1,2)},
Corollary \ref{sum} and Example \ref{1/3} imply Theorem
\ref{main}(1).

Corollary \ref{c21} and Corollary \ref{sum} imply Theorem
\ref{main}(2).
\end{proof}
\end{setup}

\begin{setup}{\bf Proof of Theorem \ref{bp}.}
\begin{proof}
Corollary \ref{sum} gives the inequalities in (1) and (2).
Examples \ref{p3}, \ref{p4} shows that both the lower bounds are
sharp.

Theorem \ref{bp} (3) is due to Corollary \ref{c22}.

Theorem \ref{bp} (4) and (5) are due to \ref{d3} through
\ref{3(1,1)}.
\end{proof}
\end{setup}

\begin{exmp}\label{p3} Fletcher found an example in \cite{C-R} with
$K^3=1$ and $p_g(X)=3$. This example of canonical 3-fold is a
hypersurface : $X_{12}\subset \bbP(1,1,1,2,6)$ with terminal
singularities. This example has 2 singularities of type
$\frac{1}{2}(1,-1,1)$. One may see that $X$ is canonically fibred
by curves of genus $2$ and $n_0(X)=6$. The 5-canonical map is
generically finite of degree $2$.
\end{exmp}

\begin{exmp}\label{p4} There is a trivial example of a 3-fold
of general type with $K^3=2$ and $p_g=4$. Take a double cover over
$\bbP^4$ with branch locus a smooth hypersurface of degree 12.
Then one gets a smooth canonical 3-fold with $K^3=2$ and $p_g=4$
and the canonical map is finite of degree 2 and $n_0(X)=5$.
\end{exmp}

\section{\bf An application to certain 4-folds}
\begin{setup}\label{fuhao} Let $Y$ be a minimal normal projective 4-fold
of general type with only canonical singularities. Assume
$p_g(Y)\ge 2$. Pick up a sub-pencil $\Lambda\subset |K_Y|$ with
$\dim \Lambda=1$. Take birational modifications $\pi:
Y'\longrightarrow Y$ like in \ref{symbols} such that the movable
part of $\pi^*(\Lambda)$ is base point free. After taking the
Stein factorization, we get an induced fibration
$f:Y'\longrightarrow B$ where $B$ is a smooth curve. Denote by
$|M|$ the movable part of $\pi^*(\Lambda)$ and by $X$ a generic
irreducible element of $|M|$. Naturally $X$ would be a smooth
projective 3-fold of general type by a weak subadditivity (see
(1.7) in \cite{Mori}).

One might hope to know whether there are parallel results to
\cite{IJM} in 4-dimensional case. This is still open because of a
pathological case. For example, take $Y$ to be a product $X\times
C$ where $X$ is a minimal 3-fold of general type with $p_g(X)=1$
and $C$ a smooth curve of genus $g\ge 2$. One would see that one
knows nothing to $Y$ since little is known to $X$. So, here, we
will exclude this case.
\end{setup}

\begin{setup}{\bf Proof of Theorem \ref{4-folds}.}
\begin{proof}
We consider the natural map
$$H^0(Y', K_{Y'})\overset\theta\longrightarrow U\subset H^0(X,
K_{Y'}|_X)\subset H^0(X, K_X)$$ where $U$ is the image of
$\theta$. If $\dim\Phi_{K_Y}(Y)\ge 2$, then $\dim_{\bC}(U)\ge 2$
and so $p_g(X)\ge 2$. If $|K_Y|$ is composed with a pencil of
3-folds $X$, by the assumption of Theorem \ref{4-folds}, we still
have $p_g(X)\ge 2$. Therefore we have $n_0(X)\le 8$ according to
\cite{IJM}.

1) Set $g:=\Phi_{\Lambda}\circ \pi$. Then $g:Y'\longrightarrow
\bbP^1$ is a surjective morphism. One sees that a general fiber of
$g$ is a smooth projective 3-dimensional scheme with each
component a 3-fold of general type. We apply Koll\'ar's technique
(\cite{Kol}) to study the canonical stability index. We have
$$\so_{\bbP^1}(1)\hookrightarrow g_*\omega_{Y'}$$
and so
$$g_*\omega_{Y'/\bbP^1}^8\hookrightarrow g_*\omega_{Y'}^{24}.$$
The standard semi-positivity says that $g_*\omega_{Y'/\bbP^1}^8$
is generated by global sections. This means that any local section
of $g_*\omega_{Y'/\bbP^1}^8$ can be glued to a global section of
$g_*\omega_{Y'}^{24}$. On the other hand, one may apply the
Kawamata-Viehweg vanishing theorem to see that $|24K_{Y'}|$ may
separate different generic irreducible elements of $|M|$. Thus the
Tankeev principle gives $n_0(Y)\le 24$.

2) We study the induced fibration $f:Y'\longrightarrow B$ which is
the Stein factorization of $g$. We may write
$$\pi^*(K_Y)\sim M+E'\equiv a_4X+E'$$
where $E'$ is an effective $\bQ$-divisor and $a_4>b=g(C)$.

The assumption implies $b=0$. We first pick up a minimal model
$X_0$ of $X$ and denote by $\sigma: X\longrightarrow X_0$ the
"contraction" onto the minimal model. Denote by $r$ the canonical
index of $X_0$. Then, according to the Base Point Free Theorem
(see Chapter 3 in \cite{KMM}), $|mrK_{X_0}|$ is base point free
for larger $m$. We have an inclusion
$$f_*\omega_{Y'/\bbP^1}^{mr}\hookrightarrow f_*\omega_{Y'}^{3mr}.$$
Since $f_*\omega_{Y'/\bbP^1}^{mr}$ is generated by global sections
and the movable part of $|mrK_{X}|$ is $mr\sigma^*(K_{X_0})$,
Lemma 2.7 of \cite{MPC} gives
$$3mr\pi^*(K_{Y})|_X\ge mr\sigma^*(K_{X_0})$$
for any larger $m$. Therefore we may find a rational number
$\beta\mapsto\frac{1}{3}$ such that $\pi^*(K_Y)|_X-\beta
\sigma^*(K_{X_0})$ is pseudo-effective. Thus we get
$$K_Y^4\ge \pi^*(K_Y)^3\cdot X\ge \frac{1}{27}K_{X_0}^3\ge
\frac{1}{81}$$ by Theorem \ref{main}. We are done.
\end{proof}
\end{setup}

\begin{rem} If there is a sound relatively minimal model program
like in 3-dimensional case, then a parallel Lemma to Lemma
\ref{replace} exists.  So in the proof above we do not need to
avoid the case $b>0$. Then Theorem \ref{4-folds} can be
strengthened to the following form:
\begin{quote}
Let $Y$ be a minimal (i.e. $K_Y$ being nef) projective 4-fold of
general type with only canonical singularities. Assume $p_g(Y)\ge
2$ and $Y$ is not canonically fibred by 3-folds of geometric genus
1. Then

(1) $n_0(Y)\le 24$;

(2) $K_Y^4\ge \frac{1}{81}$.
\end{quote}
\end{rem}
\medskip

\noindent{\bf Acknowledgment.} I would like to thank a referee so
much who gave numerous valuable suggestions and comments in order
to make the present version more readable and more well-organized.
This note was begun while I was visiting the IMS of the University
of Hong Kong in the Summer of 2004. I would like to thank both N.
Mok and X. Sun for their hospitality. Special thanks are due to
Fabrizio Catanese, Margarida Mendes Lopes and Miles Reid who
answered my frequent email queries on relevant topics. Finally I
would like to thank Y. T. Siu for suggesting this topic and
encouragements.


\end{document}